\documentclass{amsart}
\usepackage{amsfonts}
\usepackage{amsmath}
\usepackage{amssymb}
\usepackage{amsthm}

\newtheorem{lemma}{Lemma}
\newtheorem{theorem}{Theorem}
\newtheorem{rem}{Remark}
\newtheorem{exmp}{Example}
\newtheorem{cor}{Corollary}

\begin{document}

\title{On geometry of linear involutions}
\author{Mark  Pankov}
\address{Institute of Mathematics NASU,
Tereshchenkivska 3, Kiev 01601, Ukraine}
\email{pankov@imath.kiev.ua}

\subjclass[2000]{51M35, 14M15}
\keywords{Grassmann space, semilinear mapping, linear group}

\begin{abstract}
Let $V$ be an $n$-dimensional left vector space over a division ring $R$ and $n\ge 3$.
Denote by ${\mathcal G}_{k}$ the Grassmann space of $k$-dimensional subspaces of $V$
and put ${\mathfrak G}_{k}$ for the set of all pairs
$(S,U)\in {\mathcal G}_{k}\times {\mathcal G}_{n-k}$
such that $S+U=V$.
We study bijective transformations of ${\mathfrak G}_{k}$
preserving the class of base subsets and show that
these mappings are induced by semilinear isomorphisms of $V$
to itself or to the dual space $V^{*}$ if $n\ne 2k$;
for $n=2k$ this fails.
This result can be formulated as the following:
if $n\ne 2k$ and the characteristic of $R$ is not equal to
$2$ then any commutativity preserving transformation of
the set of $(k,n-k)$-involutions is extended to an automorphism
of the group {\rm GL}(V).
\end{abstract}

accepted to Advances in Geometry
\maketitle

\section{Introduction}

\subsection{}
Let $V$ be an $n$-dimensional left vector space over a division ring $R$ and $n\ge 3$.
Put ${\mathcal G}_{k}$ for the Grassmann space of $k$-dimensional subspaces of $V$.

If $B$ is a base for $V$ then
the set consisting of all $k$-dimensional subspaces spanned by vectors belonging to $B$
is called the {\it base subset} of ${\mathcal G}_{k}$ associated with $B$,
see \cite{Pankov1}, \cite{Pankov2}.
Any bijective transformation of ${\mathcal G}_{k}$ sending base subsets to
base subsets is induced by a semilinear isomorphism of $V$ to itself
or to the dual space $V^{*}$
(the second possibility can be realized only for the case when $n=2k$).
If $k=1,n-1$ then this is the Fundamental Theorem of
Projective Geometry.
For the case when $1<k<n-1$ it was established by author \cite{Pankov1},
\cite{Pankov2} (a more general result can be found in \cite{Pankov3}).

The base subsets of ${\mathcal G}_{k}$
are closely related with apartments of the Tits building associated with $V$ \cite{Tits}
(any apartment of this building consists of all flags spanned by vectors of
a certain base).
Apartment preserving transformations of the chamber sets of spherical buildings
can be extended to automorphisms of the corresponding complexes;
it follows from results of P. Abramenko and H. Van Maldeghem \cite{AVM}.

\subsection{}
Now denote by ${\mathfrak G}_{k}$ the set of all pairs
$$(S,U)\in {\mathcal G}_{k}\times {\mathcal G}_{n-k}$$
such that $S+U=V$.
Let $B$ be a base for $V$.
Consider the set of all pairs $(S,U)\in {\mathfrak G}_{k}$ such that
$S$ and $U$ are spanned by vectors belonging to the base $B$.
This set will be called the {\it base subset} of ${\mathfrak G}_{k}$
associated with $B$ (or defined by $B$).
It consists of $\binom{n}{k}$ elements and
its projections onto ${\mathcal G}_{k}$ and ${\mathcal G}_{n-k}$
are base subsets.

For any $\alpha=(S,U)\in {\mathfrak G}_{k}$ the {\it opposite element}
$\alpha^{op}=(U,S)$ belongs to ${\mathfrak G}_{n-k}$.
There is the natural bijection $p_{k}:{\mathfrak G}_{k}\to {\mathfrak G}_{n-k}$
sending each element of ${\mathfrak G}_{k}$ to the opposite element;
it transfers base subsets to base subsets.

In this paper we study bijective transformations
of ${\mathfrak G}_{k}$ preserving the class of base subsets.
The following transformations satisfy this condition:
\begin{enumerate}
\item[(1)] Any semilinear automorphism $l:V\to V$
induces the bijective transformation of ${\mathfrak G}_{k}$
sending $(S,U)$ to $(l(S),l(U))$.
\item[(2)]
There is the natural bijection of ${\mathcal G}_{i}$ onto
the Grassmann space consisting of $(n-i)$-dimensional subspaces of
the dual vector space $V^{*}$,
it maps each subspace $T$ to the annihilator $T^{0}$.
Hence any semilinear isomorphism $s:V\to V^{*}$
induces the bijective transformation of
${\mathfrak G}_{k}$ sending $(S,U)$ to $(s(U)^{0},s(S)^{0})$.
\end{enumerate}
The main result of this paper (Theorem 1) says that
{\it if $n\ne 2k$ then any bijective transformation of ${\mathfrak G}_{k}$
preserving the class of base subsets is induced by
a semilinear isomorphism of $V$ to itself or to the dual space $V^{*}$}.
If $k=1,n-1$ then this is a simple consequence of
G. W. Mackey's results \cite{Mackey} (see \cite{R} or chapter 4 of \cite{D2}).

Adjacency preserving transformations of ${\mathfrak G}_{k}$
were determined by H. Havlicek and M. Pankov \cite{HP},
but the main idea of \cite{Pankov1}, \cite{Pankov2} and \cite{Pankov3}
(characterizations of adjacency in terms of base subsets)
can not be used for the present case and
Theorem 1 will be proved by other methods.

\subsection{}
Suppose that the characteristic of $R$ is not equal to $2$.
Then for any involution $u\in {\rm GL}(V)$
there exist two  invariant subspaces $S_{+}(u)$ and $S_{-}(u)$ such that
$$u(x)=x\;\mbox{ if }\;x\in S_{+}(u)\,,
\;\;\;u(x)=-x\;\mbox{ if }\;x\in S_{-}(u)$$
and
$$V=S_{+}(u)+S_{-}(u).$$
If the dimensions of $S_{+}(u)$ and $S_{-}(u)$ are equal to $k$ and $n-k$ (respectively)
then we say that $u$ is a $(k,n-k)$-{\it involution}.
The set of all $(k,n-k)$-involutions will be denoted by ${\mathfrak I}_{k}$.
There is the natural one-to-one correspondence between
${\mathfrak I}_{k}$ and ${\mathfrak G}_{k}$.

It was mentioned in chapter 4 of \cite{D2}
(see also section 3 of \cite{D1})
that a subset of ${\mathfrak I}_{k}$ is a
maximal set of mutually permutable $(k,n-k)$-involutions
if and only if the corresponding subset of ${\mathfrak G}_{k}$ is a base subset.
This means that a bijection $f:{\mathfrak I}_{k}\to {\mathfrak I}_{k}$
is commutativity preserving
($f$ and $f^{-1}$ map commutative involutions to commutative involutions)
if and only if it can be considered as a transformation of
${\mathfrak G}_{k}$ preserving the class of base subsets.

J. Dieudonn\'{e} \cite{D1} and C. E. Rickart \cite{R} have used
Mackey's ideas \cite{Mackey} to study automorphisms of classical groups.

Our result says that {\it any commutativity preserving
bijective transformation of ${\mathfrak I}_{k}$
can be extended to an automorphism of ${\rm GL}(V)$ if $n\ne 2k$.}

\section{Results}

\begin{theorem}
If $n\ne 2k$ then any bijective transformation of ${\mathfrak G}_{k}$
preserving the class of base subsets is induced by a semilinear
isomorophism of $V$ to itself or to the dual space $V^{*}$.
\end{theorem}

If $f$ is a bijective transformation of ${\mathfrak G}_{k}$
preserving the class of base subsets
then $p_{k}fp_{n-k}$ is a bijective transformation of ${\mathfrak G}_{n-k}$
satisfying the similar condition.
Thus we need to prove Theorem 1 only for the case when $k<n-k$.
By C. E. Rickart \cite{R} (see also chapter 4 of \cite{D2}),
the required statement follows from Mackey's result \cite{Mackey} if $k=1$.
The case when $1<k<n-k$ will be considered in Section 3.

\begin{cor}
If the characteristic of $R$ is not equal to $2$ and
$n\ne 2k$ then any commutativity preserving bijection
$f:{\mathfrak I}_{k}\to {\mathfrak I}_{k}$
can be extended to an automorphism of the group ${\rm GL}(V)$.
\end{cor}

\begin{proof}[Proof of Corollary]
Let us consider $f$ as a transformation of ${\mathfrak G}_{k}$
preserving the class of base subsets.
If this mapping is induced by a semilinear automorphism
$l:V\to V$ then for any involution $u\in {\mathfrak I}_{k}$
we have
$$S_{+}(f(u))=l(S_{+}(u))\;\mbox{ and }\;S_{-}(f(u))=l(S_{-}(u))$$
and the required automorphism of ${\rm GL}(V)$ is defined by
the formula $u\to lul^{-1}$.
Now suppose that our mapping is induced by a semilinear isomorphism
$s:V\to V^{*}$.
Then
$$S_{+}(f(u))=s(S_{-}(u))^{0}\;\mbox{ and }\;S_{-}(f(u))=s(S_{+}(u))^{0}$$
and $f$ is the restriction of the automorphism $s\to s^{-1}{\check u}s$,
where ${\check u}\in {\rm GL}(V^{*})$ is the contragradient of $u$.
\end{proof}

Let $n=2k$.
Then for any element of ${\mathfrak G}_{k}$ the opposite element belongs to
${\mathfrak G}_{k}$ and Theorem 1 does not hold.
We take any subset ${\mathcal X}\subset {\mathfrak G}_{k}$
such that $\alpha \in {\mathcal X}$ implies that $\alpha^{op} \in {\mathcal X}$
(${\mathcal X}$ may be empty) and consider
the transformation of ${\mathfrak G}_{k}$ sending
each element of ${\mathcal X}$ to the opposite element and
leaving fixed elements of ${\mathfrak G}_{k}\setminus {\mathcal X}$;
if ${\mathcal X}$ coincides with ${\mathfrak G}_{k}$ then we get $p_{k}$.
Denote by $Op$ the group of all such transformations.
These transformations preserve the class of base subsets,
but non-identical elements of $Op$ are not induced by semilinear isomorphisms.

If $n=2k$ then we put $\overline{{\mathfrak G}}_{k}$
for the set of all two-element subsets $\{\alpha, \alpha^{op}\}\subset {\mathfrak G}_{k}$.
Each element of $Op$
gives the identical transformation of $\overline{{\mathfrak G}}_{k}$.
Semilinear isomorphisms of $V$ to itself or to the dual space $V^{*}$
induce bijective transformations of $\overline{{\mathfrak G}}_{k}$
and the following statement holds true.

\begin{theorem}
Let $n=2k\ge 8$ and $f$ be a bijective transformation of
${\mathfrak G}_{k}$ preserving the class of base subsets.
Then $f$ preserves the relation of being opposite
and defines a certain transformation of $\overline{{\mathfrak G}}_{k}$.
The latter mapping is induced a semilinear isomorphism of $V$ to itself or to
the dual space $V^{*}$.
\end{theorem}

In other words, if $n=2k\ge 8$ then
the group of transformations of ${\mathfrak G}_{k}$
preserving the class of base subsets is spanned by
the group $Op$ and all transformations induced by semilinear isomorphisms;
the kernel of the action of this group on $\overline{{\mathfrak G}}_{k}$
is $Op$.

\begin{cor}
Let $n=2k\ge 8$ and the characteristic of $R$ is not equal to $2$.
Let also $f$ be a commutativity preserving bijective transformation of ${\mathfrak I}_{k}$.
Then there exists an automorphism $f'$ of the group ${\rm GL}(V)$
such that $f(u)=\pm f'(u)$ for each $u\in {\mathfrak I}_{k}$.
\end{cor}
Theorem 2 is not proved for the case when $n=2k$ is equal to $4$ or $6$.

\section{Proof of Theorems 1 and 2}
Throughout the section we suppose that $1<k\le n-k$ and $n\ge 5$.
For the case when $n=2k$ we also require that $k\ge 4$.

\subsection{Main idea of the proof}
Let $\alpha=(Q,T)\in {\mathfrak G}_{m}$.
Then $(S,U)\in {\mathfrak G}_{k}$ will be called $(+)$-{\it incident}
(or $(-)$-{\it incident}) to $\alpha$ if $S$ is incident to $Q$ and $U$ is incident to $T$
(or $U$ is incident to $Q$ and $S$ is incident to $T$);
for each of these cases $\alpha$ and $(S,U)$ are said to be {\it incident}.
Put ${\mathfrak G}^{+}_{k}(\alpha)$ and ${\mathfrak G}^{-}_{k}(\alpha)$
for the sets of all elements of ${\mathfrak G}_{k}$ which are
$(+)$-incident or $(-)$-incident to $\alpha$ (respectively).
Then
$${\mathfrak G}_{k}(\alpha):={\mathfrak G}^{+}_{k}(\alpha)\cup {\mathfrak G}^{-}_{k}(\alpha)$$
consists of all elements of ${\mathfrak G}_{k}$ incident to $\alpha$.

\begin{lemma}
There exists a bijection $g:{\mathfrak G}_{k-1}\to {\mathfrak G}_{k-1}$
preserving the class of base subsets and
such that for each $\alpha\in{\mathfrak G}_{k-1}$ we have
$$f({\mathfrak G}^{+}_{k}(\alpha))={\mathfrak G}^{+}_{k}(g(\alpha))\;
\mbox{ if }\;k<n-k$$
or
$$f({\mathfrak G}_{k}(\alpha))={\mathfrak G}_{k}(g(\alpha))\;
\mbox{ if }\;n=2k.$$
\end{lemma}
This statement will be proved later.
Now we show that Theorems 1 and 2 follow from Lemma 1.

Let $\alpha=(S,U)\in{\mathfrak G}_{k}$.
Take $\beta_{i}=(Q_{i},T_{i})\in {\mathfrak G}_{k-1}$, $i=1,2$
such that
$$S=Q_{1}+Q_{2}\;\mbox{ and }\;U=T_{1}\cap T_{2}.$$
If $k<n-k$ then
$${\mathfrak G}^{+}_{k}(\beta_{1})\cap {\mathfrak G}^{+}_{k}(\beta_{2})=
\{\alpha \}.$$
Consider the case $n=2k$.
Since $Q_{1}\not\subset T_{2}$ and $Q_{2}\not\subset T_{1}$,
there are not elements of ${\mathfrak G}_{k}$ which are $(+)$-incident to
one of $\beta_{i}$ and $(-)$-incident to the other; hence
$${\mathfrak G}_{k}(\beta_{1})\cap {\mathfrak G}_{k}(\beta_{2})=
\{\alpha, \alpha^{op}\}.$$
\begin{rem}\rm{
If $n=2k=4$ then $Q_{1}, Q_{2}$ are $1$-dimensional and
the latter equality does not hold for the case when $Q_{1}\subset T_{2}$ and
$Q_{2}\subset T_{1}$ (there exist elements of ${\mathfrak G}_{k}$ which are
$(+)$-incident to $\beta_{1}$ and $(-)$-incident to $\beta_{2}$).
}\end{rem}

By Lemma 1,
\begin{equation}
{\mathfrak G}^{+}_{k}(g(\beta_{1}))\cap {\mathfrak G}^{+}_{k}(g(\beta_{2}))=
\{f(\alpha) \}\;\mbox{ if }\;k<n-k,
\end{equation}
\begin{equation}
{\mathfrak G}_{k}(g(\beta_{1}))\cap {\mathfrak G}_{k}(g(\beta_{2}))=
\{f(\alpha), f(\alpha^{op})\}\;\mbox{ if }\;n=2k.
\end{equation}
Suppose that $g$ is induced by a semilinear automorphism $l:V\to V$.
Then
$$g(\beta_{i})=(l(Q_{i}),l(T_{i}))\;\;i=1,2$$
and (1) shows that
$$f(\alpha)=(l(Q_{1})+l(Q_{2}),l(T_{1})\cap l(T_{2}))=(l(S),l(U))$$
for $k<n-k$.
If $n=2k$ then (2) implies that
$f(\alpha)$ coincides with $(l(S),l(U))$ or $(l(U),l(S))$.

Now suppose that $g$ is induced by a semilinear isomorphism $s:V\to V^{*}$.
Then
$$g(\beta_{i})=(s(T_{i})^{0},s(Q_{i})^{0})\;\;i=1,2$$
and (1) guarantees that
$$f(\alpha)=(s(T_{1})^{0}+s(T_{2})^{0},s(Q_{1})^{0}\cap s(Q_{2})^{0})=
(s(U)^{0},s(S)^{0})$$
if $k<n-k$.
By (2), $f(\alpha)$ coincides with $(s(U)^{0},s(S)^{0})$
or $(s(S)^{0},s(U)^{0})$ if $n=2k$.

Thus Theorem 1 can be proved by induction and Theorem 2 follows from Theorem 1.

\subsection{Inexact subsets of base sets}
Let ${\mathfrak B}_{k}$ be a base subset of ${\mathfrak G}_{k}$ and
$B=\{x\}^{n}_{i=1}$ be a base for $V$ associated with ${\mathfrak B}_{k}$.
Then
$$P_{1}:=Rx_{1},\dots,P_{n}:=Rx_{n}$$
is a frame.
For each $m\in \{1,\dots, n-1\}$ we denote  by ${\mathfrak B}_{m}$
the base subset of ${\mathfrak G}_{m}$ defined by  the base $B$;
in what follows such base subsets will be called {\it associated}
with ${\mathfrak B}_{k}$.

If $\alpha\in {\mathfrak B}_{m}$ then we put
${\mathfrak B}_{k}(\alpha)$, ${\mathfrak B}^{+}_{k}(\alpha)$ and
${\mathfrak B}^{-}_{k}(\alpha)$
for the intersections of ${\mathfrak B}_{k}$ with ${\mathfrak G}_{k}(\alpha)$,
${\mathfrak G}^{+}_{k}(\alpha)$ and ${\mathfrak G}^{-}_{k}(\alpha)$ (respectively).

\begin{rem}{\rm
Let $(S,U)\in {\mathfrak B}_{i}$ and $(Q,T)\in {\mathfrak B}_{j}$.
If one of the subspaces $S$ or $U$ is incident to $Q$ or $T$ then
the other subspace is incident $T$ or $Q$ (respectively) and
the pairs $(S,U)$ and $(Q,T)$ are incident.
}\end{rem}

A subset ${\mathfrak R}\subset {\mathfrak B}_{k}$ is said
to be {\it exact} if it is contained in exactly one base subset of ${\mathfrak G}_{k}$;
otherwise, we will say that ${\mathfrak R}$ is {\it inexact}.
Let ${\mathcal R}_{k}$ and ${\mathcal R}_{n-k}$
be the projections of ${\mathfrak R}$ onto
the Grassmann spaces ${\mathcal G}_{k}$ and ${\mathcal G}_{n-k}$.
For any $i\in \{1,\dots, n\}$ we denote by $S_{i}({\mathfrak R})$ the intersection
of all elements of ${\mathcal R}_{k}\cup {\mathcal R}_{n-k}$
containing $P_{i}$.
If $(S,U)\in {\mathfrak B}_{k}$ then one of the subspaces $S$ or $U$ contains $P_{i}$;
this means that each $S_{i}({\mathfrak R})$ is non-zero if ${\mathfrak R}$ is not empty.
It is trivial that ${\mathfrak R}$ is exact if and only if each
$S_{i}({\mathfrak R})$ coincides with $P_{i}$.

\begin{exmp}{\rm
Let $\beta=(Q,T)\in {\mathfrak B}_{2}$ and $Q=P_{i}+P_{j}$.
Then $S_{p}({\mathfrak B}_{k}(\beta))=P_{p}$ for all $p\ne i,j$
and the subspaces $S_{i}({\mathfrak B}_{k}(\beta))$ and $S_{j}({\mathfrak B}_{k}(\beta))$
are coincident with $Q$.
Thus ${\mathfrak B}_{k}(\beta)$ is inexact;
moreover, it is a maximal inexact subset of ${\mathfrak B}_{k}$:
for any $(S,U)\in {\mathfrak B}_{k}-{\mathfrak B}_{k}(\beta)$
the subset
$$\{(S,U)\}\cup {\mathfrak B}_{k}(\beta)$$
is exact
(indeed, $S$ intersects $Q$ only
by one of the subspaces $P_{i}$, $P_{j}$ and $U$ intersects $Q$ by the other).
}\end{exmp}

\begin{lemma}
If ${\mathfrak R}$ is a maximal inexact subset of ${\mathfrak B}_{k}$
then there exists $\beta\in {\mathfrak B}_{2}$ such that
${\mathfrak R}={\mathfrak B}_{k}(\beta)$.
\end{lemma}

\begin{proof}
Since ${\mathfrak R}$ is inexact, for some number $i$
the dimension of $S_{i}({\mathfrak R})$ is not less than $2$.
Denote this dimension by $m$.
There exists a unique $(n-m)$-dimensional subspace $T$
such that $(S_{i}({\mathfrak R}), T)$ is a element of ${\mathfrak B}_{m}$.
We define
$$\gamma:=
\begin{cases}
(S_{i}({\mathfrak R}), T)\mbox{ if } m\le n-m\\
(T, S_{i}({\mathfrak R}))\mbox{ if } m>n-m.
\end{cases}
$$
and get ${\mathfrak R}\subset {\mathfrak B}_{k}(\gamma)$.
Take any $\beta\in {\mathfrak B}_{2}$ which is $(+)$-incident to $\gamma$.
Then ${\mathfrak R}\subset{\mathfrak B}_{k}(\beta)$.
Since our inexact set is maximal, we have the inverse inclusion.
\end{proof}

\subsection{}
Let $U$ and $U'$ be $m$-dimensional subspaces of $V$, $1<m<n-1$.
Recall that the {\it distance} $d(U,U')$ between $U$ and $U'$ is equal to
$$m-\dim U\cap U'=\dim (U+U')-m;$$
if $d(U,U')=1$ then $U$ and $U'$ are said to be {\it adjacent}.
The distance $d(U,U')$ is the minimal number $i$ such that there exists a sequence of
$m$-dimensional subspaces
$$U=U_{0},U_{1},\dots,U_{i}=U'$$
where $U_{j-1}$ and $U_{j}$ are adjacent for each $j\in \{1,\dots,i\}$.

For any $\alpha=(S,U)$ and $\beta=(S',U')$ belonging  to ${\mathfrak B}_{k}$
$$d(S,S')=d(U,U')$$
is not greater than $k$ and $(\alpha,\beta)$ will be called an
$i$-{\it pair} if this number is equal to $i$.
For this case we will also say that
the {\it distance} $d(\alpha,\beta)$ is equal to $i$.
It is easy to see that $d(\alpha,\beta)$
is the minimal number $i$ such that there is a sequence
$$\alpha=\alpha_{0},\alpha_{1},\dots,\alpha_{i}=\beta$$
of elements of ${\mathfrak B}_{k}$ where $\alpha_{j-1}$ and $\alpha_{j}$
form a $1$-pair for each $j\in \{1,\dots,i\}$.

By our hypothesis,
${\mathfrak B}'_{k}:=f({\mathfrak B}_{k})$ is a base subset.
For each number $m\in \{1,\dots, n-1\}$
we denote by ${\mathfrak B}'_{m}$ the base subsets of ${\mathfrak G}_{m}$ associated
with ${\mathfrak B}'_{k}$.

Since ${\mathfrak R}\subset {\mathfrak B}_{k}$ is a maximal inexact subset
if and only if $f({\mathfrak R})$ is a maximal inexact subset of ${\mathfrak B}'_{k}$,
Lemma 2 implies that
for each $\alpha\in {\mathfrak B}_{2}$
there exists $\alpha'\in {\mathfrak B}'_{2}$ such that
$$f({\mathfrak B}_{k}(\alpha))={\mathfrak B}'_{k}(\alpha');$$
there is only one $\alpha'$ satisfying this equality
(indeed, the condition $n\ge 5$ guarantees that
for distinct $\alpha, \beta \in {\mathfrak B}_{2}$ the sets
${\mathfrak B}_{k}(\alpha)$ and ${\mathfrak B}_{k}(\beta)$ are not coincident).
We set $f_{2}(\alpha):=\alpha'$, then
$f_{2}:{\mathfrak B}_{2}\to {\mathfrak B}'_{2}$ is a bijection
(if $n=2k=4$ then for any $\alpha\in {\mathfrak B}_{2}$
the set ${\mathfrak B}_{k}(\alpha)$ coincides with ${\mathfrak B}_{k}(\alpha^{op})$
and the mapping $f_{2}$ is not well defined).

\begin{lemma}
The bijection $f_{2}$ preserves the distance.
\end{lemma}

\begin{proof}
Let $\alpha_{1}=(Q_{1},T_{1})$ and $\alpha_{2} =(Q_{2},T_{2})$
be distinct elements of ${\mathfrak B}_{2}$ and $i$ be the distance between them.
Then $i=1$ ($\alpha_{1}$ and $\alpha_{2}$ form a $1$-pair)
or $i=2$ ($\alpha_{1}$ and $\alpha_{2}$ are incident);
for each of these cases, the cardinal number of
the intersection of ${\mathfrak B}_{k}(\alpha_{1})$ and ${\mathfrak B}_{k}(\alpha_{2})$
will be denoted by $c_{1}$ and $c_{2}$, respectively.

Let $i=1$. Then
$$\beta:=(Q_{1}+Q_{2},T_{1}\cap T_{2})$$
is an element of ${\mathfrak B}_{3}$ and
$(S,U)\in {\mathfrak B}_{k}$ belongs to
${\mathfrak B}_{k}(\alpha_{1})\cap {\mathfrak B}_{k}(\alpha_{2})$
if and only if one of the following conditions is fulfilled:
\begin{enumerate}
\item[(a)] $Q_{1}+Q_{2}\subset S$,
\item[(b)] $S\subset T_{1}\cap T_{2}$
\end{enumerate}
(see Remark 2).
There are exactly $\binom{n-3}{k-3}$
and $\binom{n-3}{k}$ distinct $(S,U)\in {\mathfrak B}_{k}$
satisfying (a) and (b), respectively.
Therefore,
$$c_{1}=\binom{n-3}{k-3}+\binom{n-3}{k}.$$
Now suppose that $i=2$.
Then $Q_{1}\subset T_{2}$, $Q_{2}\subset T_{1}$ and
$\beta$ is element of ${\mathfrak B}_{4}$ (recall that $n\ge 5$).
It is easy to see that $(S,U)\in {\mathfrak B}_{k}$ belongs to
${\mathfrak B}_{k}(\alpha_{1})\cap {\mathfrak B}_{k}(\alpha_{2})$
if and only if one of the following conditions holds:
\begin{enumerate}
\item[(a)] $Q_{1}+Q_{2}\subset S$,
\item[(b)] $S\subset T_{1}\cap T_{2}$,
\item[(c)] $S$ is incident to both $Q_{1}$ and $T_{2}$,
\item[(d)] $S$ is incident to both $Q_{2}$ and $T_{1}$.
\end{enumerate}
There are $\binom{n-4}{k-4}$, $\binom{n-4}{k}$,
$\binom{n-4}{k-2}$, $\binom{n-4}{k-2}$
distinct $(S,U)\in {\mathfrak B}_{k}$ satisfying (a), (b), (c), (d), respectively
(if $k=2,3$ then some of these numbers are equal to $0$).
This implies that
$$c_{2}=\binom{n-4}{k-4}+\binom{n-4}{k}+2\binom{n-4}{k-2}.$$
We have
$$c_{1}=\binom{n-3}{k-3}+\binom{n-3}{k}=
\binom{n-4}{k-4}+\binom{n-4}{k-3}+\binom{n-4}{k-1}+\binom{n-4}{k}$$
and
$$c_{2}-c_{1}=2\binom{n-4}{k-2}-\binom{n-4}{k-3}-\binom{n-4}{k-1}>0.$$
Since
$$f({\mathfrak B}_{k}(\alpha_{1})\cap {\mathfrak B}_{k}(\alpha_{2}))=
{\mathfrak B}'_{k}(f_{2}(\alpha_{1}))\cap {\mathfrak B}'_{k}(f_{2}(\alpha_{2})),$$
the sets ${\mathfrak B}_{k}(\alpha_{1})\cap {\mathfrak B}_{k}(\alpha_{2})$
and ${\mathfrak B}'_{k}(f_{2}(\alpha_{1}))\cap {\mathfrak B}'_{k}(f_{2}(\alpha_{2}))$
have the same cardinal number which is equal to $c_{1}$ or $c_{2}$.
We have $c_{2}\ne c_{1}$, this means that $f_{2}$ preserves the distance.
\end{proof}

\begin{lemma}
Let $\alpha=(S,U)\in {\mathfrak B}_{m}$, $2\le m\le k$.
Let also
\begin{equation}
\beta_{1}=(Q_{1},T_{1}),\dots,\beta_{m-1}=(Q_{m-1},T_{m-1})
\end{equation}
be a sequence of elements of ${\mathfrak B}_{2}$ such that
$\beta_{i}$ and $\beta_{j}$ form a $1$-pair if and only if $j=i\pm 1$,
$$
S=Q_{1}+\dots +Q_{m-1}\;\mbox{ and }\; U=T_{1}\cap\dots\cap T_{m-1}.
$$
Then
$${\mathfrak B}_{k}(\alpha)=
{\mathfrak B}_{k}(\beta_{1})\cap\dots\cap {\mathfrak B}_{k}(\beta_{m-1}).$$
\end{lemma}

\begin{proof}
Simple induction.
\end{proof}

\begin{lemma}
If $2\le m \le k$ then for each $\alpha\in {\mathfrak B}_{m}$
there exists $\alpha'\in {\mathfrak B}'_{m}$ such that
$$f({\mathfrak B}_{k}(\alpha))={\mathfrak B}'_{k}(\alpha').$$
\end{lemma}

\begin{proof}
For any $\alpha=(S,U)\in {\mathfrak B}_{m}$ there is a sequence (3) satisfying
the conditions of the previous lemma.
By Lemma 3, the analogous assertion holds for the sequence
$$f_{2}(\beta_{1})=(Q'_{1},T'_{1}),\dots,f_{2}(\beta_{m-1})=(Q'_{k-1},T'_{m-1});$$
i.e. $f_{2}(\beta_{i})$ and $f_{2}(\beta_{j})$ form a $1$-pair if and only if $j=i\pm 1$.
Then
$$\alpha':=(Q'_{1}+\dots +Q'_{m-1}, T'_{1}\cap\dots\cap T'_{m-1})$$
is an element of ${\mathfrak B}'_{m}$ and  Lemma 4 shows that
$$f({\mathfrak B}_{k}(\alpha))=
f({\mathfrak B}_{k}(\beta_{1})\cap\dots\cap {\mathfrak B}_{k}(\beta_{m-1}))=$$
$${\mathfrak B}'_{k}(f_{2}(\beta_{1}))\cap\dots\cap {\mathfrak B}'_{k}(f_{2}(\beta_{m-1}))=
{\mathfrak B}'_{k}(\alpha').$$
\end{proof}

Now suppose that $\alpha \in {\mathfrak B}_{k}$.
If $k<n-k$ then there is unique $\alpha'\in {\mathfrak B}'_{k}$
satisfying the equality of Lemma 5 and
we define $f_{k}(\alpha):=\alpha'$.
The mapping $f_{k}:{\mathfrak B}_{k}\to {\mathfrak B}'_{k}$
is bijective
(if $n=2k$ then ${\mathfrak B}_{k}(\alpha)={\mathfrak B}_{k}(\alpha^{op})$
for any $\alpha \in {\mathfrak B}_{k}$ and $f_{k}$ is not well defined).

\begin{lemma}
If $k<n-k$ then $f_{k}$ preserves the distance.
\end{lemma}

\begin{proof}
Let $\alpha=(S,U)$ and $\beta=(Q,T)$ be distinct elements of ${\mathfrak B}_{k}$ and
$i$ be the distance between them.
For each of the cases $i=1,\dots, k$ the cardinal number of the set
${\mathfrak B}_{k}(\alpha)\cap {\mathfrak B}_{k}(\beta)$
will be denoted by $c_{1},\dots, c_{k}$ (respectively).
If $i<k$ then ${\mathfrak B}_{k}(\alpha)\cap {\mathfrak B}_{k}(\beta)$
consists only of all $(M,N)\in {\mathfrak B}_{k}$ such that $M$
is contained in $U\cap T$.
The dimension of $U\cap T$ is equal to $n-k-i$ and
$$c_{i}=
\binom{n-k-i}{k}\;\mbox{ if }\;i<k.$$
If $i=k$ then $\alpha$ and $\beta$ are incident;
hence they are both belonging to ${\mathfrak B}_{k}(\alpha)\cap {\mathfrak B}_{k}(\beta)$;
moreover, $(M,N)\in {\mathfrak B}_{k}$ is an element of
${\mathfrak B}_{k}(\alpha)\cap {\mathfrak B}_{k}(\beta)$ if
$M$ is contained in $U\cap T$.
The subspace $U\cap T$ is $(n-2k)$-dimensional and
$$c_{k}=\binom{n-2k}{k}+2.$$
A direct verification shows that $c_{1}>0$ is not equal to $c_{2},\dots, c_{k}$.
Since
$$f({\mathfrak B}_{k}(\alpha)\cap {\mathfrak B}_{k}(\beta))=
{\mathfrak B}'_{k}(f_{k}(\alpha))\cap {\mathfrak B}'_{k}(f_{k}(\beta)),$$
the sets ${\mathfrak B}_{k}(\alpha)\cap {\mathfrak B}_{k}(\beta)$
and ${\mathfrak B}'_{k}(f_{k}(\alpha))\cap {\mathfrak B}'_{k}(f_{k}(\beta))$
have the same cardinal number.
This means that $f_{k}$ preserves the class of $1$-pairs and the claim follows.
\end{proof}

\begin{lemma}
If $k<n-k$ then $f_{k}$ is the restriction of $f$ to ${\mathfrak B}_{k}$ and
Lemma {\rm 6} implies that two elements of ${\mathfrak B}_{k}$ form a $1$-pair
if and only if their $f$-images form a $1$-pair.
\end{lemma}

\begin{proof}
Let $\alpha\in {\mathfrak B}_{k}$.
First of all note that
$$f_{k}({\mathfrak B}_{k}(\alpha))={\mathfrak B}'_{k}(f(\alpha)).$$
Indeed, we have the following chain of equivalences
$$\beta\in {\mathfrak B}_{k}(\alpha)\Leftrightarrow \alpha\in {\mathfrak B}_{k}(\beta)
\Leftrightarrow f(\alpha)\in {\mathfrak B}'_{k}(f_{k}(\beta)) \Leftrightarrow
f_{k}(\beta)\in {\mathfrak B}'_{k}(f(\alpha)).$$
On the other hand, the set ${\mathfrak B}_{k}(\alpha)$ consists
of all $\beta \in {\mathfrak B}_{k}$ incident to $\alpha$.
Since $f_{k}$ is distance preserving
and two distinct elements of ${\mathfrak B}_{k}$
are incident if and only if the distance between them is $k$,
the set $f_{k}({\mathfrak B}_{k}(\alpha))$ consists of all
elements of ${\mathfrak B}'_{k}$ incident to $f_{k}(\alpha)$.
This implies the equality
$$f_{k}({\mathfrak B}_{k}(\alpha))={\mathfrak B}'_{k}(f_{k}(\alpha))$$
and we get $f(\alpha)=f_{k}(\alpha)$.
\end{proof}

\begin{lemma}
If $k<n-k$ then for any $\alpha\in{\mathfrak B}_{k-1}$
there exists $\alpha'\in{\mathfrak B}_{k-1}$
$$f({\mathfrak B}^{+}_{k}(\alpha))={\mathfrak B}'^{+}_{k}(\alpha').$$
\end{lemma}

\begin{proof}
Let $\alpha=(Q,T)$ be an element of ${\mathfrak B}_{k-1}$.
We say that ${\mathfrak R}\subset {\mathfrak B}_{k}(\alpha)$
is a {\it $1$-subset} of ${\mathfrak B}_{k}(\alpha)$
if any two distinct elements of ${\mathfrak R}$ form a $1$-pair
and ${\mathfrak R}$ is a maximal subset of ${\mathfrak B}_{k}(\alpha)$
satisfying this condition.
Clearly, ${\mathfrak B}^{+}_{k}(\alpha)$ is a $1$-subset.

There are exactly $p:=\binom{n-k+1}{k+1}$ distinct $(k+1)$-dimensional
subspaces of $V$ which are spanned by vectors of the base $B$
and contained in $T$ (since $k<n-k$, we have $p>1$).
Denote these subspaces by $M_{1},\dots,M_{p}$
and put ${\mathfrak R}_{i}$ for the set of all
$(S,U)\in {\mathfrak B}_{k}(\alpha)$ such that $S\subset M_{i}$.
Then each ${\mathfrak R}_{i}$ is a $1$-subset.
Any $1$-subset of ${\mathfrak B}_{k}(\alpha)$
coincides with ${\mathfrak B}^{+}_{k}(\alpha)$ or certain ${\mathfrak R}_{i}$.
The set ${\mathfrak B}^{+}_{k}(\alpha)$ can be characterized
by the following property:
it is the unique $1$-subset having empty intersection with all other $1$-subsets
of ${\mathfrak B}_{k}(\alpha)$
(indeed, if $M_{i}\cap M_{j}$ is $k$-dimensional
then ${\mathfrak R}_{i}\cap {\mathfrak R}_{j}\ne \emptyset$).

By Lemma 5, we have
$$f({\mathfrak B}_{k}(\alpha))={\mathfrak B}'_{k}(\alpha')$$
for certain $\alpha'\in{\mathfrak B}_{k-1}$.
Lemma 7 shows that $f$ transfers $1$-subsets of ${\mathfrak B}_{k}(\alpha)$ to
$1$-subsets of ${\mathfrak B}'_{k}(\alpha')$.
Then $f({\mathfrak B}^{+}_{k}(\alpha))$
coincides with ${\mathfrak B}'^{+}_{k}(\alpha')$.
\end{proof}

\subsection{Proof of Lemma 1 for the case $k<n-k$}

Let $k<n-k$ and $\alpha=(Q,T)\in {\mathfrak G}_{k-1}$.
We say that ${\mathfrak B}\subset {\mathfrak G}^{+}_{k}(\alpha)$
is a {\it base subset} of ${\mathfrak G}^{+}_{k}(\alpha)$ if
there exists a base subset ${\mathfrak B}_{k}$ of ${\mathfrak G}_{k}$
such that $\alpha$ belongs to the associated base subset of
${\mathfrak G}_{k-1}$ and ${\mathfrak B}={\mathfrak B}^{+}_{k}(\alpha)$.

Any base subset
$$\{(S_{1},U_{1}),\dots,(S_{n-k+1},U_{n-k+1})\}$$
of ${\mathfrak G}^{+}_{k}(\alpha)$
defines the frame $T\cap S_{1},\dots, T\cap S_{n-k+1}$ for the subspace $T$.
Conversely, any frame ${\mathcal F}=\{P_{1},\dots,P_{n-k+1}\}$ for $T$
gives the base subset
$${\mathfrak B}_{\mathcal F}:=\{(Q+P_{1},P_{2}+\dots+P_{n-k+1}),\dots,
(Q+P_{n-k+1},P_{1}+\dots+P_{n-k})\}$$
of ${\mathfrak G}^{+}_{k}(\alpha)$.

Let ${\mathfrak B}$ and  ${\mathfrak B}'$ be base subsets of ${\mathfrak G}^{+}_{k}(\alpha)$.
We write ${\mathfrak B}\sim {\mathfrak B}'$ if
the intersection of ${\mathfrak B}$ and ${\mathfrak B}'$
contains at least $2$ elements.
For the case when there exists a sequence of base subsets
$${\mathfrak B}={\mathfrak B}_{1}\sim{\mathfrak B}_{2}\sim\dots\sim
{\mathfrak B}_{i}={\mathfrak B}'$$
of ${\mathfrak G}^{+}_{k}(\alpha)$ we write
${\mathfrak B}\simeq {\mathfrak B}'$.

\begin{lemma}
Let ${\mathfrak B}$ and  ${\mathfrak B}'$ be base subsets of ${\mathfrak G}^{+}_{k}(\alpha)$.
Lemma 8 implies the existence of
$\beta,\beta'\in {\mathfrak G}_{k-1}$ such that
$$f({\mathfrak B})\subset{\mathfrak G}^{+}_{k}(\beta)\;\mbox{ and }\;
f({\mathfrak B}')\subset{\mathfrak G}^{+}_{k}(\beta').$$
If ${\mathfrak B}\simeq {\mathfrak B}'$ then $\beta=\beta'$.
\end{lemma}

\begin{proof}
It is clear that we can restrict ourself only to the case when
${\mathfrak B}\sim {\mathfrak B}'$.
If $(S_{1},U_{1})$ and $(S_{2},U_{2})$ are the $f$-images of
two distinct elements of ${\mathfrak B}\cap{\mathfrak B}'$
then $\beta=(S_{1}\cap S_{2}, U_{1}+U_{2})=\beta'$.
\end{proof}

\begin{lemma}
${\mathfrak B}\simeq {\mathfrak B}'$ for any two base subsets
${\mathfrak B}$ and ${\mathfrak B}'$ of ${\mathfrak G}^{+}_{k}(\alpha)$.
\end{lemma}

\begin{proof}
Let ${\mathcal F}$ and ${\mathcal F}'$ be the frames for $T$
associated with ${\mathfrak B}$ and ${\mathfrak B}'$
(${\mathfrak B}={\mathfrak B}_{{\mathcal F}}$ and
${\mathfrak B}'={\mathfrak B}_{{\mathcal F}'}$).
First we consider the case when
$${\mathcal F}=\{P_{1},\dots,P_{n-k},P_{n-k+1}\}\;\mbox{ and }\;
{\mathcal F}'=\{P_{1},\dots,P_{n-k},P'\}.$$
If $P'\subset P_{n-k}+P_{n-k+1}$ then the condition $n-k\ge 3$
guarantees that
$$(Q+P_{1}, P_{2}+P_{3}+\dots+P_{n-k}+P_{n-k+1})=(Q+P_{1}, P_{2}+P_{3}+\dots+P_{n-k}+P'),$$
$$(Q+P_{2}, P_{1}+P_{3}+\dots+P_{n-k}+P_{n-k+1})=(Q+P_{2}, P_{1}+P_{3}+\dots+P_{n-k}+P')$$
belong to ${\mathfrak B}\cap{\mathfrak B}'$ and
${\mathfrak B}\sim {\mathfrak B}'$.

Now suppose that we have ${\mathfrak B}\simeq {\mathfrak B}'$ if
$P'\subset P_{i+1}+\dots +P_{n-k+1}$ and consider the case when
$P'\subset P_{i}+\dots +P_{n-k+1}$.
The subspace $P_{i}+P'$ intersects $P_{i+1}+\dots +P_{n-k+1}$ by
a $1$-dimensional subspace $P''$ such that
$P_{1},\dots, P_{n-k}, P''$ is a frame for $T$.
Let ${\mathfrak B}''$ be the base subset of
${\mathfrak G}^{+}_{k}(\alpha)$ associated with the latter frame.
Then ${\mathfrak B}\simeq {\mathfrak B}''$, ${\mathfrak B}'\sim {\mathfrak B}''$
and we get ${\mathfrak B}\simeq {\mathfrak B}'$.

If $|{\mathcal F}\cap {\mathcal F}'|<n-k$ then
there exists a sequences of frames
$${\mathcal F}={\mathcal F}_{0}, {\mathcal F}_{1},\dots,{\mathcal F}_{i}={\mathcal F}'$$
for $T$ such that $|{\mathcal F}_{j-1}\cap {\mathcal F}_{j}|=n-k$
for each $j\in \{1,\dots,i\}$
(the number $i$ is equal to $n-k+1-|{\mathcal F}\cap {\mathcal F}'|$),
we have
$${\mathfrak B}={\mathfrak B}_{{\mathcal F}_{0}}\simeq
{\mathfrak B}_{{\mathcal F}_{1}}\simeq\dots\simeq
{\mathfrak B}_{{\mathcal F}_{i}}={\mathfrak B}'$$
and ${\mathfrak B}\simeq {\mathfrak B}'$.
\end{proof}

Each element of ${\mathfrak G}^{+}_{k}(\alpha)$ is contained in a base
subset of ${\mathfrak G}^{+}_{k}(\alpha)$ and Lemmas 9, 10 imply the existence
of $\beta\in {\mathfrak G}_{k-1}$ such that
$$f({\mathfrak G}^{+}_{k}(\alpha))\subset {\mathfrak G}^{+}_{k}(\beta).$$
Since $f^{-1}$ preserves the class of base subsets of ${\mathfrak G}_{k}$,
we have the inverse inclusion.

We define $g(\alpha):=\beta$. The mapping $g$ is bijective.
If ${\mathfrak B}_{k-1}$ is a base subset of ${\mathfrak G}_{k-1}$
and ${\mathfrak B}_{k}$ is the base subset of ${\mathfrak G}_{k}$ associated with
${\mathfrak B}_{k-1}$ then $g({\mathfrak B}_{k-1})$ is the base subset of ${\mathfrak G}_{k-1}$
associated with $f({\mathfrak B}_{k})$.
Similar arguments show that $g^{-1}$ maps base subsets to base subsets.

\subsection{Proof of Lemma 1 for the case $n=2k$}

Now suppose that $n=2k\ge 8$ and  $\alpha=(Q,T)\in {\mathfrak G}_{k-1}$.
We say that ${\mathfrak B}\subset {\mathfrak G}_{k}(\alpha)$
is a {\it base subset} of ${\mathfrak G}_{k}(\alpha)$ if
there exists a base subset ${\mathfrak B}_{k}$ of ${\mathfrak G}_{k}$
such that $\alpha$ belongs to the associated base subset of
${\mathfrak G}_{k-1}$ and ${\mathfrak B}={\mathfrak B}_{k}(\alpha)$.
There is a one-to-one correspondence between
base subsets of ${\mathfrak G}_{k}(\alpha)$ and frames for $T$.

Let ${\mathfrak B}$ and  ${\mathfrak B}'$ be base subsets of ${\mathfrak G}_{k}(\alpha)$.
We write ${\mathfrak B}\sim {\mathfrak B}'$ if
the intersection of ${\mathfrak B}$ and ${\mathfrak B}'$
contains at least $6$  elements
(if $(S,U)$ belongs to a base subset of ${\mathfrak G}_{k}(\alpha)$
then the same holds for $(U,S)$).
If there exists a sequence of base subsets
$${\mathfrak B}={\mathfrak B}_{1}\sim{\mathfrak B}_{2}\sim\dots\sim
{\mathfrak B}_{i}={\mathfrak B}'$$
of ${\mathfrak G}_{k}(\alpha)$ then we write
${\mathfrak B}\simeq {\mathfrak B}'$.

\begin{lemma}
Let ${\mathfrak B}$ and  ${\mathfrak B}'$ be base subsets of ${\mathfrak G}_{k}(\alpha)$.
By Lemma 5, there exist
$\beta,\beta'\in {\mathfrak G}_{k-1}$ such that
$$f({\mathfrak B})\subset{\mathfrak G}_{k}(\beta)\;\mbox{ and }\;
f({\mathfrak B}')\subset{\mathfrak G}_{k}(\beta').$$
If ${\mathfrak B}\simeq {\mathfrak B}'$ then $\beta=\beta'$.
\end{lemma}

\begin{proof}
We need to prove this equality only for the case when ${\mathfrak B}\sim {\mathfrak B}'$.
Suppose that $\beta=(M,N)$.
We choose three elements $(S_{i},U_{i})$, $i=1,2,3$
of $f({\mathfrak B})\cap f({\mathfrak B}')$ such that
$$M=S_{1}\cap S_{2}\cap S_{3}$$
and $N$ is the sum of any $U_{p}$ and $U_{q}$ if $p\ne q$.
Since $f({\mathfrak B})$ is a base subset of ${\mathfrak G}_{k}(\beta)$
(Lemma 5),
$$\dim U_{1}\cap U_{2}\cap U_{3}=k-2$$
and any $S_{p}$ and $U_{q}$ are not adjacent.
This means that $\beta'=(M,N)$.
\end{proof}

\begin{rem}{\rm
The assertion of Lemma 11 can not be proved if
${\mathfrak B}\cap {\mathfrak B}'$ has only $4$ elements and
${\mathfrak B}\not\simeq{\mathfrak B}'$.
Indeed, let $\beta=(M,N)\in {\mathfrak G}_{k-1}$
and $(S_{i},U_{i})$, $i=1,2$ be elements of ${\mathfrak G}_{k}(\beta)$
such that
$$M=S_{1}\cap S_{2}\;\mbox{ and }\; N=U_{1}+U_{2};$$
then
$$\beta':=(U_{1}\cap U_{2},S_{1}+S_{2})$$
is another element of ${\mathfrak G}_{k-1}$ and
$(S_{i},U_{i})$, $i=1,2$ belong to ${\mathfrak G}_{k}(\beta')$.
}\end{rem}

\begin{lemma}
${\mathfrak B}\simeq {\mathfrak B}'$ for any two base subsets
${\mathfrak B}$ and ${\mathfrak B}'$ of ${\mathfrak G}_{k}(\alpha)$.
\end{lemma}

\begin{proof}
Let ${\mathcal F}$ and ${\mathcal F}'$ be the frames for $T$
associated with ${\mathfrak B}$ and ${\mathfrak B}'$.
Suppose that
$${\mathcal F}=\{P_{1},\dots,P_{n-k},P_{n-k+1}\},\;\;
{\mathcal F}'=\{P_{1},\dots,P_{n-k},P'\}$$
and $P'\subset P_{n-k}+P_{n-k+1}$.
Since $n-k\ge 4$,
$$(Q+P_{1}, P_{2}+P_{3}+\dots+P_{n-k}+P_{n-k+1})=(Q+P_{1}, P_{2}+P_{3}+\dots+P_{n-k}+P'),$$
$$(Q+P_{2}, P_{1}+P_{3}+\dots+P_{n-k}+P_{n-k+1})=(Q+P_{2}, P_{1}+P_{3}+\dots+P_{n-k}+P'),$$
$$(Q+P_{3}, P_{1}+P_{2}+P_{4}+\dots+P_{n-k}+P_{n-k+1})=
(Q+P_{3}, P_{1}+P_{2}+P_{4}+\dots+P_{n-k}+P')$$
and the opposite elements belong to ${\mathfrak B}\cap{\mathfrak B}'$ and
${\mathfrak B}\sim {\mathfrak B}'$.
For other cases the proof is similar to the proof of Lemma 10.
\end{proof}

Lemmas 11 and 12 together with arguments of the previous subsection give the claim.

\begin{rem}{\rm
If $n=2k$ is equal to $6$ or $4$
then the intersection of two distinct base subsets of
${\mathfrak G}_{k}(\alpha)$ contains at most two non-opposite elements
($k=3$) or two elements which are opposite ($k=2$).
Theorem 2 is not proved for these cases.
}\end{rem}

\vspace{3ex}
\noindent
{\bf Acknowledgement}.
I thank Krzysztof Pra{\.z}mowski and Mariusz {\.Z}ynel
for useful discussions during my staying in University of Bia{\l}ystok.


\begin{thebibliography}{99}

\bibitem{AVM}
P. Abramenko  and H. Van Maldeghem,
On oppositions in spherical buildings and twin buildings.
Ann. Combinatorics 4(2000), 125 -- 137.

\bibitem{D1}
J. Dieudonn\'{e}, On the automorphisms of the classical groups.
Memoirs Amer. Math. Soc. 2(1951), 1 -- 95.

\bibitem{D2}
J. Dieudonn\'{e},
La G\'{e}om\'{e}trie des Groupes Classiques.
Springer-Verlag, Berlin, 1971.

\bibitem{Mackey} G. W. Mackey,
Isomorphisms of  normed linear spaces. Ann. of Math. 43 (1942),
244 -- 260.

\bibitem{HP}
H. Havlicek  and M. Pankov,
Transformations on the product of Grassmann spaces.
Demonstratio Math. XXXVIII (2005), to appear.

\bibitem{Pankov1} M. Pankov,
Transformations of Grassmannians and automorphisms of linear groups.
J. Geom. 75(2002), 132 -- 150.

\bibitem{Pankov2} M. Pankov,
Transformations of Grassmannians
preserving the class of base sets.
J. Geom. 79(2004), 169 -- 176.

\bibitem{Pankov3} M. Pankov,
A characterization of geometrical mappings of Grassmann spaces,
Result. Math. 45(2004), 319 -- 327.


\bibitem{R} C. E. Rickart,
Isomorphic groups of linear transformations,
Amer. J. Math. 72(1950), 451 -- 464.

\bibitem{Tits} J. Tits,
Buildings of spherical type and finite BN-pairs, Lect. Notes in Math.
386, Springer-Verlag, 1974.
\end{thebibliography}
\end{document}